\documentclass[12pt]{article}
\usepackage{fourier}
\usepackage{index}
\usepackage[expansion=false]{microtype}
\usepackage[reqno]{amsmath}
\usepackage{fixmath}
\usepackage{amssymb, amsthm, enumerate}
\usepackage{mathtools}
\usepackage{pgf,tikz,pgfplots}
\usetikzlibrary{automata, positioning, arrows}
\usepackage[abbrev,msc-links]{amsrefs} 
\usepackage{subcaption}
\usepackage{diagbox}

\usepackage{authblk}

\usepackage{hyperref}

\newtheoremstyle{plainsl}%
	{\topsep}
	{\topsep}
	{\slshape} 
	{}
	{\normalfont\bfseries}
	{.}
	{ }
	{}

\swapnumbers

{\theoremstyle{plainsl}
\newtheorem{theorem}{Theorem}[section]
\newtheorem*{theorem*}{Theorem}
\newtheorem{lemma}[theorem]{Lemma}
\newtheorem{corollary}[theorem]{Corollary}}
{\theoremstyle{remark}
\newtheorem{example}[theorem]{Example}}
\newtheorem{remark}[theorem]{Remark}

\newcommand\cref[1]{Corollary~\ref{cor:#1}}

\renewcommand\proof{\noindent\textsl{Proof. }}
\newcommand\sqr[2]{{\vbox{\hrule height.#2pt
    \hbox{\vrule width.#2pt height#1pt \kern#1pt
        \vrule width.#2pt}\hrule height.#2pt}}}
\renewcommand\qed{%
	\ifmmode\eqno\sqr53
	\else\nolinebreak\ \hfill\sqr53\medbreak\fi}

\numberwithin{equation}{section}

\newcommand\comp[1]{{\mkern2mu\overline{\mkern-2mu#1}}}

\newcommand\seq[3]{#1_{#2},\ldots,#1_{#3}}

\DeclareMathOperator{\rk}{rk}
\DeclareMathOperator{\tr}{tr}
\DeclareMathOperator{\col}{col}

\DeclareMathOperator{\diag}{diag}
\DeclareMathOperator{\adj}{adj}
\DeclareMathOperator{\summ}{sum}

\newcommand\ip[2]{\langle #1,#2\rangle}

\begin{document}
\title{Cospectral graphs obtained by edge deletion}
\author[1]{Chris Godsil\footnote{Research supported by 
Natural Sciences and Engineering Research Council of Canada, 
Grant No. RGPIN-9439}}
\author[1,2]{Wanting Sun
\footnote{Visiting Ph. D. student at 
the Department of Combinatorics \& Optimization at University of Waterloo 
from March 2022 to March 2023; 
supported by China Scholarship Council (No. 202106770031)}}
\author[1]{Xiaohong Zhang}

\affil[1]{Department of Combinatorics \& Optimization, 
University of Waterloo, Waterloo, Ontario, Canada}
\affil[2]{Faculty of Mathematics and Statistics, Central China Normal University, Wuhan, 430079, PR China}

\maketitle

\begin{abstract}
We show how to construct graphs which are simultaneously cospectral with respect to several symmetric matrices that encode adjacency. Our main operation consists of locating cliques in highly regular graphs (for example, strongly regular graphs) and removing edges in distinct ways preserving the spectrum of the adjacency, Laplacian, unsigned Laplacian, and normalized Laplacian matrix. We also provide  a more general construction of Laplacian cospectral graphs. This has applications in Laplacian state transfer in the study of quantum information transfer. 
\end{abstract}

\section{Introduction}

There are a number of useful matrices we can associate with a graph $X$. We have, of course,
the adjacency matrix $A=A(X)$. If $D=D(X)$ is the diagonal matrix with $D_{i,i}$ equal to the
valency of the vertex $i$, we also have the \textsl{Laplacian} $L(X)=D-A$, the 
\textsl{unsigned Laplacian} (also called \textsl{signless Laplacian}) $S(X)=D+A$,  and the \textsl{normalized Laplacian} 
$N(X)=D^{-1/2}L(X)D^{-1/2}$, where $D^{-1}$ is the  Moore–Penrose inverse of $D$. 
It has long been known that, although it is a useful invariant, the spectrum of $A(X)$ does not determine $X$. In this paper we provide new constructions of families of graphs where any two graphs in a family 
 are cospectral with respect to all the four matrices.

Our construction makes use of a type of highly regular graphs, 1-walk regular graphs. A graph is \textsl{walk-regular}  
if for any positive integer $k$, 
the number of closed walks of length $k$ is the same at all vertices. 
If further, 
the number of walks from vertex $u$ to $v$ of length $k$ is the same for all adjacent vertex pairs $u, v$, 
then we say $X$ is \textsl{1-walk regular}. 
For example, strongly regular graphs are 1-walk regular. 
In the next section, we will discuss 1-walk regular graphs in more details. 
Given a graph $X =(V,E)$, if $F\subseteq E$, 
we use $X\backslash F$ to denote the graph obtained from $X$ by deleting all the edges in $F$. 
One of our main results is the following (Theorem~\ref{thm:cosALQsubg}).

\begin{theorem*}
	Let $X_1$ and $X_2$ be two cospectral 1-walk regular graphs. 
For $i=1,2$,  let $Y_i$ be a graph contained in a clique of $X_i$. Assume $Y_1\cong Y_2$.  
Then $X_1\backslash E(Y_1)$ and $X_2\backslash E(Y_2)$, and their complements, 
are cospectral
with respect to the adjacency, the Laplacian,  the unsigned Laplacian, and the normalized Laplacian matrices. 
Furthermore, $X_1\backslash E(Y_1)$ and $X_2\backslash E(Y_2)$ are not isomorphic if $X_1$ are $X_2$ are not.  \qed
\end{theorem*}

For a concrete example, there exists a strongly regular graph $X$ with 150 edges and parameters $(25, 12, 5, 6)$, such that 
all the 150 graphs obtained by deleting an edge from $X$
are pairwise non-isomorphic and pairwise cospectral with respect to the adjacency matrix $A$, the Laplacian matrix $L$, the unsigned Laplacian $S$, and the normalized Laplacian $N$ (here $X_1=X_2=X$). The complements of these 150 graphs are also cospectral with respect to these four matrices.

We discuss some related earlier results. 
In \cite{GodsilMcKay1982}, Godsil and McKay gave several constructions for adjacency cospectral graphs. 
 Dalf\'o, van Dam and Fiol \cite{DvF2012} constructed adjacency cospectral graphs by considering perturbations of almost distance-regular graphs. 
The perturbations include vertex deletion, edge deletion, amalgamating vertices, and adding pendant vertices. 
Some of our results on adjacency cospectrality of graphs overlap with the results by Dalf\'o et al.~\cite{DvF2012}.  
Dutta constructed unsigned Laplacian cospectral graphs in \cite{unsignedLcos2020} and
Butler \cite{ButleradjnorL} constructed graphs that are cospectral with respect to adjacency matrix and normalized Laplacian matrix. In \cite[Corollary 1.1]{butler2022complements}, Butler et al.~constructed graphs that are cospectral with respect to adjacency matrix, Laplacian, unsigned and normalized Laplacian matrices: Let $H_1$ and $H_2$ be two regular cospectral graphs, then for any rooted graph $G$, the rooted products $H_1\circ G$ and $H_2\circ G$ are cospectral with respect to the four matrices. In \cite[Theorem~4.1]{WangLiLu}, Wang et al.~gave a construction of pairs of graphs cospectral
with respect to all the four matrices,  by using a special case of Godsil-McKay switching.

Let $\Delta(X)$ denote the distance matrix of a graph $X$, and let $\comp{X}$ denote the complement of $X$. In 1977, McKay \cite{McKaycostree} constructed  pairs of non-isomorphic trees that are not distinguished by the spectrum of $A(T)$, $A(\comp{T})$, $L(T)$, $S(T)$ and $\Delta(T)$. (His full list of matrices is longer.) Osborne \cite[Theorem~3.3.2]{Osborne} constructed pairs of trees that are not distinguished by the spectrum of $A(T)$, $A(\comp{T})$, $L(T)$, $S(T)$, and $N(T)$. 
  
The graphs we construct are related by edge-deletion. For vertex deletion, we have the
following. Let $X$ be a graph with adjacency matrix $A$.  
For any subset $R$ of $V(X)$, 
denote the induced subgraph of $X$ on $V(X)\backslash R$ by $X\backslash R$, and the 
characteristic polynomial of $X$ with respect to the adjacency matrix as $\phi(X,t)$. 
Jacobi's Theorem (Theorem 1.2 from Chapter~4 of \cite{GodsilAC}) connects characteristic polynomials of graphs and determinant of a principal submatrix of the inverse of a matrix: 
\begin{equation}\label{eq:Jacobithm}
	\det \left(\big((tI-A)^{-1}\big)_{R,R}\right)=\frac{\phi(X\backslash R,t)}{\phi(X,t)}. 
\end{equation}
This implies that if $X$ is a strongly regular graph, 
and $Y_1$, $Y_2$ are induced subgraphs of $X$ that are (adjacency) cospectral with cospectral complements, 
then so are $X\backslash Y_1$ and $X\backslash Y_2$.

Our above theorem says that 
if $Y_1$ and $Y_2$ are two isomorphic graphs that both lie in a clique of a 1-walk regular graph $X$ (the cliques can be different ), then the two graphs obtained from $X$ by removing the edges of $Y_1$ or edges of $Y_2$, respectively, 
are cospectral with respect to all the four matrices $A$, $L$, $S$, and $N$. 
Relaxing the condition by only requiring $Y_1$ and $Y_2$ to be Laplacian cospectral (not necessarily isomorphic), we further show that the resulting graphs are Laplacian cospectral. 
Similarly, assume $Y_1$ and $Y_2$, and their complements, are cospectral with respect to the unsigned Laplacian matrix.   
Then the two graphs obtained by removing the edges of $Y_1$ or edges of $Y_2$, respectively, from a clique of a 1-walk regular graph, 
are cospectral with respect to the unsigned Laplacian matrix.

This work is motivated by quantum information transfer.  
A quantum spin system of $n$ interacting qubits can be modelled by a graph on $n$ vertices. 
The state of the quantum system is represented by a density matrix $\mathcal{D}$  (positive semidefinite of trace 1) of size $n\times n$. 
Under the Heisenberg model, the Hamiltonian of the system is $L(X)$, and 
the transition matrix of the system at time $t$ is the unitary matrix $U(t)=e^{itL(X)}$. Given two states $\mathcal{D}_1$ and $\mathcal{D}_2$, if there is time $t$ such that $U(t)\mathcal{D}_1U(-t) =\mathcal{D}_2$, then we say there is \textsl{perfect state transfer} from $\mathcal{D}_1$ to  $\mathcal{D}_2$. 
Connections between the graph structures of $X$ and possible state transfer has been studied \cite{godsil2017realstatetransfer}. 
Denote $e_1,\ldots, e_n$ the standard basis vectors. 
Interesting states include: vertex state $\mathcal{D}=e_1e_1^T$ \cite{whenpstGodsil}, pair state $\mathcal{D}=\frac12 (e_1-e_2)(e_1-e_2)^T$ \cite{pair}, and Laplacian state $\mathcal{D}_Y=\frac{1}{2m}L(Y)$, where $Y$ is some graph on $n$ vertices with $m$ edges \cite{Lapstate}. 
In particular, 
if $Y_1$ and $Y_2$ are two spanning subgraphs of $X$ with the same number of edges,  there is perfect \textsl{Laplacian state} transfer from $Y_1$ to $Y_2$ if   
 there is a time $t$ such that $e^{itL(X)}L(Y_1)e^{-itL(X)}=L(Y_2)$. 
Since $U(t)=e^{itL(X)}$ and $L(X)$ commute, 
we have $U(t)L(X)U(-t) =L(X)$, 
and hence $U(t)\big(L(X)-L(Y_1)\big)U(-t) =L(X)-L(Y_2)$. 
Therefore if Laplacian state transfer occurs between $Y_1$ and $Y_2$, 
then in addition to $L(Y_1)$ and $L(Y_2)$ being similar, 
$L(X)-L(Y_1)$ and $L(X)-L(Y_2)$ are also similar. 
That is, the edge-deleted subgraphs $X\backslash E(Y_1)$ and $X\backslash E(Y_2)$ are Laplacian cospectral. 
One can similarly use the unsigned Laplacian as the Hamiltonian, and a scaled unsigned Laplacian matrix as a state.

\section{1-walk regular graphs}

Recall that a graph is 1-walk regular
if it is walk-regular and for any positive integer $k$, 
the number of walks from vertex $u$ to $v$ of length $k$ is the same for all adjacent vertex pairs $u, v$. 
Let $A$ and  $B$ be two matrices of the same size, say $m\times n$, 
then their \textsl{Schur product} $A\circ B$ is the matrix of size $m\times n$ such that 
\[
(A\circ B)_{j,k}=A_{j,k}B_{j,k}. 
\]
In terms of this matrix product, 
a graph $X$ with adjacency matrix $A$ is 1-walk regular if and only if, 
 for any positive integer $k$, 
 there exist scalars $a_k, b_k$ such that 
\[
A^k \circ I = a_k I \text{ and } A^k \circ A = b_k A.
\]
Distance-regular graphs are 1-walk regular; 
more generally, 
given a symmetric association scheme $\mathcal{A}=\{A_0=I, A_1,\ldots, A_d\}$, the graph with adjacency matrix $A_i$ is  1-walk regular for $i=1,\ldots, d$. 
More details about association schemes can be found  \cite[Chapter~12]{GodsilAC}.

Assume $M$ is a Hermitian matrix with exactly $m$ distinct eigenvalues $\theta_1,\ldots, \theta_m$.  
For $r=1,\ldots, m$, 
let $E_r$ denote the orthogonal projection matrix onto the $\theta_r$-eigenspace of $M$. 
Then 
\[
M =\sum_{r=1}^m \theta_r E_r. 
\]
and the above equation is called the \textsl{spectral decomposition} of $M$. 
Denote the Kronecker delta function by $\delta_{i,j}$. These projections  satisfy $E_rE_s=\delta_{r,s}E_r$. 
As a consequence, 
if $f(x)$ is a function which is defined on each eigenvalue of $M$,
then we may assume that $f(M)=\sum_r f(\theta_r)E_r$ (our usage follows Chapter 9 of \cite{LancasterTismenetskyM}).  
In particular, 
\[
(tI-M)^{-1}=\sum_r \frac{1}{t-\theta_r} E_r.
\]

Since each of the orthogonal projection matrices $E_1,\ldots, E_m$ is a polynomial in $M$, the matrix algebra $\langle M \rangle$ generated by $M$ is the same as the one generated by $\{E_1,\ldots, E_m\}$. 
Therefore a graph $X$ with adjacency matrix $A=\sum_r \theta_r E_r$ is 1-walk regular if and only if for each $r$, 
there exist scalars $\alpha_r$ and $\beta_r$ such that 
\[
E_r \circ I = \alpha_r I, \;\; E_r\circ A = \beta_r A.
\] 

We have seen two commonly used characterizations of 1-walk regular graphs: in terms of powers of adjacency matrix, or in terms of  eigenspace projection matrices. 
Now we give more characterizations. 
In particular, 1-walk regularity is equivalent to cospectrality of certain subgraphs of $X$. 

We denote the adjugate of a matrix $M$ by $\adj(M)$.   
We use $X\backslash u$ to denote the graph obtained from a graph $X$ by deleting a vertex $u$, and $X\backslash \{u,v\}$ to denote the graph obtained from a graph $X$ by deleting vertices $u$ and $v$.

\begin{theorem}\label{thm:1walkrechar}
Let $X$ be a graph on $n$ vertices with adjacency matrix $A$. 
Assume $A=\sum_{r=1}^m \theta_r E_r$ is the spectral decomposition of $A$. 
The following are equivalent.
\begin{itemize}
\item[(a)]
$X$ is 1-walk regular.
\item[(b)]
For any positive integer $k$, 
there exist scalars $a_k$ and $b_k$ such that 
$A^k \circ I = a_k I \text{ and } A^k \circ A = b_k A$.
\item[(c)]
For any integer $k=1,\ldots, n-1$, 
there exist scalars $a_k$ and $b_k$ such that 
$A^k \circ I = a_k I \text{ and } A^k \circ A = b_k A$.
\item[(d)]
For $r=1,\ldots, m$, 
there exist scalars $\alpha_r$ and $\beta_r$ such that 
$E_r \circ I = \alpha_r I$ and $E_r\circ A = \beta_r A$. 
\item[(e)]
There exist functions $a(t)$ and $b(t)$ such that 
$(tI-A)^{-1} \circ I = a(t) I \text{ and } (tI-A)^{-1} \circ A = b(t) A$.
\item[(f)]
There exist polynomials $f(t)$ and $g(t)$ such that 
$\adj(tI - A) \circ I = f(t) I$ and  $\adj(tI-A) \circ A = g(t) A$. 
\item[(g)]
$\phi(X\backslash u,t)$ is the same for all $u\in V(X)$, 
and $\phi(X\backslash \{u,v\},t)$ is the same for adjacent vertices $u,v$ in $X$.  
\end{itemize}
\end{theorem}

\proof
The equivalence between $(a), (b)$ and $(d)$ has been shown. 
The equivalence of $(b)$ and $(c)$ follows from Cayley-Hamilton Theorem. 
The equivalence of $(d)$ and $(e)$ follows from 
$(tI-A)^{-1} =\sum_{r=1}^m \frac{1}{t-\theta_r}E_r$.
The equivalence of $(e)$ and $(f)$ follows from 
$ \adj(tI-A) = \phi(X,t) (tI-A)^{-1}$.

Finally we prove the equivalence of $(f)$ and $(g)$. 
Note that $[\adj(tI-A)]_{u,u}= \phi(X\backslash u,t)$,
the equivalence of 
$\phi(X\backslash u,t)$ is the same for all $u\in V(X)$ 
and $\adj(tI - A) \circ I = f(t) I$ for some $f(t)$ follows. 
Denote $[\adj(tI-A)]_{u,v}$ as $\phi_{uv}(X,t)$, 
by Corollary~1.4 from Chapter 4 of \cite{GodsilAC}, 
$\phi_{uv}(X,t)$ can be expressed 
in terms of the characteristic polynomials of subgraphs of $X$ with at most two vertices deleted. 
In fact, 
for any two vertices $u,v$ of $X$, 
\begin{equation}\label{eq:phiuv}
\phi_{u,v}(X,t)=\sqrt{\phi(X\backslash u,t) \phi(X\backslash v,t) - 
\phi(X,t) \phi(X\backslash \{u,v\},t)}.
\end{equation}
Therefore for a walk-regular graph $X$, 
$\phi(X\backslash \{u,v\},t)$ is the same for all adjacent vertex pairs $u$ and $v$ if and only if $\phi_{u,v}(X,t)$ is, 
and hence $(f)$ and $(g)$ are equivalent. \qed

 Note the equivalence of condition $(g)$ to 1-walk regularity was also obtained by Dalf\'o, van Dam and Fiol under the notion of 1-cospectrality and 1-isospectrality  \cite{DvF2012}, which is also a consequence of Proposition 6 of Godsil \cite{Godsilequi}.

Let $X$ be a 1-walk regular graph and let $e$ be an edge in $X$. 
The scalars in Theorem~\ref{thm:1walkrechar} play an important role in determining the characteristic polynomial of the edge-deleted subgraph $X\backslash e$ with respect to $A,L,S$ or $N$. 
It turns out cospectral 1-walk regular graphs share these parameters. 
Since 1-walk regular graphs are regular, cospectrality with respect to any of these four matrices are equivalent.  

\begin{lemma}\label{lem:cosp1walkreg}
Let $X$ and $X'$ be two 1-walk regular graphs. 
If $X$ and $X'$ are cospectral, then they share the same parameters in Theorem~\ref{thm:1walkrechar}. 
\end{lemma}
\proof 
Assume $X$ has $n$ vertices. 
Denote the adjacency matrices of $X$ and $X'$ by $A$ and $A'$, respectively. 
If $A$ and $A'$ are cospectral with eigenvalues $\lambda_1,\ldots, \lambda_n$, then for any positive integer $k$, 
$\tr(A^k) =\sum_{i=1}^n \lambda_i^k = \tr((A')^k)$. 
Since $A^k$ and $(A')^k$ have constant diagonals, we know that the scalars $a_k$'s in Theorem~\ref{thm:1walkrechar} $(b)$ are shared by $X$ and $X'$. \\
Denote the sum of all the entries of a matrix $B$ by $\summ(B)$. 
Then 
\[
\summ(A\circ A^k)=\tr(A^{k+1}) =\tr((A')^{k+1})=\summ(A'\circ (A')^k). 
\]
Since $A\circ A^k = b_k A$ for some constant $b_k$, we know that the scalars $b_k$'s in Theorem~\ref{thm:1walkrechar} $(b)$ are shared by $X$ and $X'$ as well. It is straightforward that if $X$ and $X'$ share one set of parameters from  Theorem~\ref{thm:1walkrechar}, say the $a_k$'s and $b_k$'s, then they also share other sets of parameters in the theorem. \qed

Let $X$ and $X'$ be two cospectral 1-walk regular graphs. 
By Lemma~\ref{lem:cosp1walkreg}, if $u,v$ are two adjacent vertices of $X$ and $u',v'$ are two adjacent vertices of $X'$, then $\phi_{u,v}(X,t)=\phi_{u',v'}(X',t)$. 
By equation \eqref{eq:phiuv}, we know the corresponding one-vertex or adjacent two-vertex deleted subgraphs of $X$ and $X'$ are cospectral as well.  
\begin{corollary}\label{coro:1walkresamepara}
Let $X$ and $X'$ be two cospectral 1-walk 
regular graphs. 
For any edge $uv$ of $X$ and edge $u'v'$ of $X'$, we have 
 $\phi(X\backslash u)=\phi(X'\backslash u')$ and 
$\phi(X\backslash \{u,v\})=\phi(X'\backslash \{u',v'\})$.  \qed
\end{corollary} 

Even though the complement of a 1-walk regular graph may not be 1-walk regular, it is walk regular and has some other nice properties. 

\begin{lemma} \label{1walkregcomple}
Let $X$ be a 1-walk regular graph with adjacency matrix $A$ and complement $\bar X$. 
Denote the adjacency matrix of $\bar X$ as $\bar A$. 
Then 
\begin{itemize}
\item[(a)]
For any positive integer $k$, 
there exist scalars $a_k$ and $b_k$ such that 
$\bar A^k \circ I = a_k I \text{ and } \bar A^k \circ A = b_k A$.
\item [(b)]
Assume $\bar A=\sum_{r=1}^s \lambda_r F_r$ is the spectral decomposition of $\bar A$, 
then for $r=1,\ldots, s$, 
there exist scalars $\alpha_r$ and $\beta_r$ such that 
$F_r \circ I = \alpha_r I$ and $F_r\circ A = \beta_r A$.
\item [(c)]
There exist functions $a(t)$ and $b(t)$ such that 
$(tI-\bar A)^{-1} \circ I = a(t) I \text{ and } (tI-\bar A)^{-1} \circ A = b(t) A$.
\item[(d)]
$\phi(\bar X\backslash u,t)$ is the same for all $u\in V(X)$, 
and $\phi(\bar X \backslash \{u,v\},t)$ is the same for any edge $uv$ of $X$.  
\end{itemize}
Furthermore, if $X$ and $X'$ are cospectral 1-walk regular graphs, 
then $\bar X$ and $\bar{X'}$ share the above parameters. 
\end{lemma}

\proof 
Assume $X$ is $k$-regular, then $AJ=JA=kJ$ for the all-ones matrix $J$. 
Therefore for some scalars $c_{pq}$ and $c_p$, 
\[
\bar A^k =(J-I-A)^k=\sum_{\substack{0\leq p,q \\  p+q\leq k }}
c_{pq} A^pJ^q = \sum_{0\leq p \leq k}c_p A^p
\]
and the first claim follows from the 1-walk regularity of $X$. 
Conditions $(b)$, $(c)$ and $(d)$ are all equivalent to $(a)$, and can be proved similarly as in Theorem~\ref{thm:1walkrechar}.   \qed

Theorem~\ref{thm:1walkrechar} implies that if $X$ is 1-walk regular, 
then for a function $f(x)$ defined on the eigenvalues of the adjacency matrix $A$, 
the matrix $f(A)$ has constant diagonal, 
and is constant on the entries corresponding to edges of $X$. 
In fact a similar result hold for the Laplacian, unsigned Laplacian, and normalized Laplacian matrix. Furthermore, a similar result holds for $\bar X$. 
Recall that these above four matrices are denoted by $A(X)$, $L(X)$, $S(X)$, and $N(X)$, respectively. When the graph $X$ is clear, we write simply $A, L,S,$ and $N$.

\begin{lemma}\label{lem:1walkinver}
Let $X$ be a 1-walk regular graph. 
Let $H$ denote its adjacency, Laplacian, unsigned Laplacian, 
or normalized Laplacian matrix, and define $\bar H$ similarly for $\bar X$.  
If $f(x)$ is a function that is defined on all eigenvalues of $H$, 
then there exist scalars $\alpha_f$ and $\beta_f$ such that 
\[
f(H)\circ I = \alpha_f I \text{ and } f(H)\circ A=\beta_f A. 
\]
If $g(x)$ is a function that is defined on all eigenvalues of $\bar H$, 
then there exist scalars $\alpha_g$ and $\beta_g$ such that 
\[
f(\bar H)\circ I = \alpha_g I \text{ and } f(\bar H)\circ A=\beta_g A. 
\]
Furthermore, if $X'$ is a 1-walk regular graph that is cospectral to $X$, then the scalars $\alpha_f$, $\beta_f$, $\alpha_g$, and $\beta_g$ are the same for $X'$. 
\end{lemma} 
\proof
The result for the adjacency matrix follows from Theorem~\ref{thm:1walkrechar}. 
The other cases follows from the fact that a function in $H$ is also a function in $A$.  
In fact, assume $X$ is $d$-regular, 
then $f(L)=f(dI-A)$, 
$f(S)=f(dI+A)$, 
and $f(N)=f(I-\frac1d A)$. 
Making use of Lemma~\ref{1walkregcomple}, the claim about $\bar H$ can be proved similarly. The result about $X'$ follows by further applying Lemma~\ref{lem:cosp1walkreg}. 
\qed

\section{Constructing graphs cospectral with respect to $A$, $L$, $S$ and $N$}\label{sec:samesubgraph}

For an edge $e$ of $X$,
let $X\backslash e$ denote the graph obtained by deleting the edge $e$ from $X$. 
Let $X$ be a 1-walk regular graph. 
We know from Theorem~\ref{thm:1walkrechar} that 
all the graphs in the set $\{X\backslash u \; |\; u\in V(X)\}$ are adjacency cospectral, 
and so are all the graphs in $\{X\backslash \{u,v\} \; |\; uv\in E(X)\}$. 
We will see that deleting any edge of $X$ also results in adjacency cospectral graphs. That is, all the graphs in $\{X\backslash uv \; |\; uv\in E(X)\}$ are cospectral. 
The characteristic polynomials of these vertex-deleted or edge-deleted subgraphs of a graph are closely related.

\begin{lemma} \cite[Chapter~4, Lemma~1.5]{GodsilAC}  \label{lem:phiXdeledge}
Let $X$ be a graph and let $e=uv$ be an edge of $X$. 
Then 
\begin{align}\label{eq:phideled}
\phi(X,t)= & \phi(X\backslash e,t) - \phi(X\backslash \{u,v\},t) \\ \nonumber
& - 2\sqrt{\phi(X\backslash u,t) \phi(X\backslash v,t) - 
\phi(X\backslash e,t) \phi(X\backslash \{u,v\},t)}.
\end{align}
\end{lemma}

Now we show that different ways of removing an edge from a 1-walk regular graph result in adjacency cospectral graphs. 
Note that part of the result is a consequence of results 
 in \cite{Godsilequi, DvF2012}. 
 If $u$ and $v$ are two non-adjacent vertices of a graph $X$, 
we denote the graph obtained from $X$ by adding an edge between $u$ and $v$ as $X+uv$. 
 
\begin{theorem}\label{thm:edgeApoly}
Let $X$ and $X'$ be two cospectral 1-walk regular graphs. 
Then for any edge $e$ of $X$ and any edge $e'$ of $X'$, 
the two graphs $X\backslash e$ and $X'\backslash e'$, and their complements, are adjacency cospectral. 
\end{theorem}

\proof 
Assume $X$ has $n$ vertices and  $e=uv$. 
Solving $\phi(X\backslash e,t)$ from \eqref{eq:phideled}, 
we have
\begin{align*}
\phi(X\backslash e,t) 
& = \phi (X,t) - \phi(X\backslash \{u,v\},t)
\pm 2\sqrt{\phi(X\backslash u,t) \phi(X\backslash v,t) - 
\phi(X,t) \phi(X\backslash \{u,v\},t)}\\
&= \phi (X,t) - \phi(X\backslash \{u,v\},t) \pm 2\phi_{u,v}(X,t) 
\quad \text{ by } \eqref{eq:phiuv}. 
\end{align*}
Recall that for a graph $Y$ on $n$ vertices, 
the coefficient of $x^{n-2}$ in $\phi(Y,t)$ is equal to $-|E(Y)|$. 
Now comparing the coefficient of $x^{n-2}$ on both sides of the above equation, 
we know that only the plus sign is valid, 
that is, 
\[
\phi(X\backslash e,t) = \phi (X,t) - \phi(X\backslash \{u,v\},t) + 2\phi_{u,v}(X,t). 
\]
Now the cospectrality of  $X\backslash e$ and $X'\backslash e'$ follows from the cospectrality of $X$ and $X'$,  Lemma~\ref{lem:cosp1walkreg},  and Corollary~\ref{coro:1walkresamepara} (recall $\phi_{u,v}(X,t)=[\adj(tI-A)]_{u,v}$). 
Let $Y$ be a graph, $u,v$ be two non-adjacent vertices of $Y$, and let $e=uv$. Then 
\[
\phi(Y+ e,t) = \phi (Y,t) - \phi(Y\backslash \{u,v\},t) - 2\phi_{u,v}(Y,t). 
\]
Let $Y=\bar X$, then the claim about cospectral complements holds by Lemma~\ref{1walkregcomple}. 
\qed

\begin{example}\label{exm:deledge}
As mentioned in Section~\ref{sec:sgr25,12,5,6}, 
there are 15 non-isomorphic strongly regular graphs with parameters 
SRG(25,12, 5, 6). 
By taking $X=X'$ in Theorem~\ref{thm:edgeApoly} to be  
one of these graphs $X_i$ and by removing an edge from $X$ in different ways,  
we obtain a family of graphs 
 that are pairwise adjacency cospectral  with cospectral complements. 
In particular, 
two of the 15 strongly regular graphs have the property that 
different ways of removing an edge from the graph result in non-isomorphic graphs. 
Therefore each of the two graphs provides a family of 150 graphs such that any two graphs in the same family are cospectral with cospectral complement but not isomorphic.  
These two graphs correspond to $X_1$ and $X_3$ in Table~\ref{tab:dataek3p3k4}, with graph6-string being, respectively,
\medskip\\
'X\textasciitilde{}zfCqTc\{YPT`fUQidaeNRKxItIMpholosZFKjXHZGnDZDYHwuF',\\
 \smallskip
'X\textasciitilde{}zfCqTc\{YPR`jUQidaeNRLXIrIMphoxKsVXKixPZCnD[fBHuQl'. \\ 
Furthermore, one can combine the two families (or more generally, the 15 families obtained from the 15 strongly regular graphs) of graphs to form a bigger family of pairwise adjacency cospectral graphs (with cospectral complement), as strongly regular graphs with the same parameters are cospectral and Theorem~\ref{thm:edgeApoly} applies.  
\end{example}

In fact, as we will see in the following, 
a more general cospectral property holds when edges are removed from cospectral 1-walk regular graphs.

\subsection{Deleting subgraphs in cliques}
 
We have established that by removing an edge from cosptreal 1-walk regular graphs, 
we get a set of pairwise adjacency cospectral graphs. In fact, 
these graphs are also cospectral with respect to 
the Laplacian matrix, unsigned Laplacian matrix, 
and normalized Laplacian matrix.
We show a more general result: 
removing edges of a small graph from a clique of cospectral 1-walk regular graphs result in graphs cospectral with respect to $A, L, S$ and $N$.  

We make use of the following result about the inverse of a rank-1 update of
an invertible matrix. 

\begin{theorem}[Sherman-Morrison\cite{ShermanMorrison}]\label{thm:shermanmorrison}
Suppose $B$ is an $n\times n$ invertible real matrix 
and $u,v\in \mathbb{R}^n$.
Then $B+uv^T$ is invertible if and only if 
$1+v^TB^{-1}u\neq 0$. 
In this case,
\begin{equation}\label{eq:rank1inver}
(B+uv^T)^{-1}=B^{-1}-\frac{B^{-1}uv^TB^{-1}}{1+v^TB^{-1}u}. 
\end{equation}
\end{theorem}

For two matrices $C$ and $D$ such that both $CD$ and $DC$ are defined,  
 $\det(I-CD)$ and $\det(I-DC)$ are closely related. 
 
\begin{lemma}\label{lem:i-cddcddet}\cite{alggraph}
Assume $C$ and $D^T$ are both matrices of size $m\times n$, 
then 
\[
\det(I_m-CD)=\det(I_n-DC). 
\]
In particular, 
if $C=u$ and $D=v^T$ for some real vectors $u,v\in\mathbb{R}^n$, 
we have 
\begin{equation} \label{eq:i+uv^T}
\det(I_n-uv^T)=(1-v^Tu). 
\end{equation}
\end{lemma}

First we prove a result concerning entries of a matrix related to edge deletion from a clique of a 1-walk regular graph. We use $\bar A$ to denote $A(\bar X)$ and  $Z\cong Z'$ to indicate the two graphs $Z$ and $Z'$ are isomorphic.  

\begin{lemma}\label{lem:valueinclique}
Let $X$ and $X'$ be two cospectral 1-walk regular graphs. 
Let $u_1,v_1,\ldots, u_r,v_r$ be vertices in a clique of $X$ and let $Z$ be the (multi)graph on these vertices with edge set $\{u_1v_1,\ldots, u_rv_r\}$.  
Let $u'_1,v'_1\ldots, u'_r,v'_r$ be vertices in a clique of $X'$ and define $Z'$ similarly. Let $H$ (resp. $H'$) be the adjacency, Laplacian, unsigned Laplacian, or normalized Laplacian matrix of $X$ (resp. $X'$). Assume $Z\cong Z'$ with the map sending $u_i$ to $u'_i$ and $v_i$ to $v'_i$ for $i=1,\ldots, r$ an isomorphism between them. 
Then for any function $f(x)$ defined on the eigenvalues of $H$, 
{\small
\begin{equation*}\label{eq:edgesplus}
e_{u_{r}}^T 
\Bigl(tI-f(H)\pm e_{u_1}e_{v_1}^T \pm\cdots \pm
e_{u_{r-1}}e_{v_{r-1}}^T\Bigl)^{-1}
e_{v_{r}}
=
e_{u'_{r}}^T 
\Bigl(tI-f(H')\pm e_{u'_1}e_{v'_1}^T \pm\cdots \pm
e_{u'_{r-1}}e_{v'_{r-1}}^T\Bigl)^{-1}
e_{v'_{r}}, 
\end{equation*}
}
and for any function $\bar f(x)$ defined on the eigenvalues of $\bar H$, 
{\small
\begin{equation*}
e_{u_{r}}^T 
\Bigl(tI-\bar f(\bar H)\pm e_{u_1}e_{v_1}^T \pm\cdots \pm
e_{u_{r-1}}e_{v_{r-1}}^T\Bigl)^{-1}
e_{v_{r}}
=
e_{u'_{r}}^T 
\Bigl(tI-\bar f(\bar H')\pm e_{u'_1}e_{v'_1}^T \pm\cdots \pm
e_{u'_{r-1}}e_{v'_{r-1}}^T\Bigl)^{-1}
e_{v'_{r}}. 
\end{equation*}
}
\end{lemma}

\proof 
The function $g(x)=(t-f(x))^{-1}$ is also defined on all eigenvalues of $H$. 
By Lemma~\ref{lem:1walkinver},  there exist scalars $\alpha(t),\beta(t)$ 
such that $\left(tI-f(H)\right)^{-1} \circ I = \alpha(t) I$ and 
$\left(tI-f(H)\right)^{-1} \circ A = \beta(t) A$. 
By Lemma~\ref{lem:1walkinver}, if $X$ and $X'$ are cospectral, 
the two equations still hold for 
the same $\alpha(t)$ and $\beta(t)$ if we replace $H$ with $H'$.

We prove the result by induction. 
When $r=1$,  since $u_1$ and $v_1$ are in the same clique of $X$ and $\delta_{u_1,v_1}=\delta_{u'_1,v'_1}$, we have 
\[
e_{u_1}^T \left(tI-f(H)\right)^{-1} e_{v_1}
= \delta_{u_1,v_1} \alpha(t) + (1- \delta_{u_1,v_1}) \beta(t)
=e_{u'_1}^T \left(tI-f(H')\right)^{-1} e_{v'_1}, 
\]
and the value depends only on whether $u_1$ and $v_1$ are the same or not. 
  
For $s=1,\ldots, r$, let $Z_s$ be the (multi)graph with vertex set $V(X)$ and edge set $\{u_1v_1, \ldots, u_sv_s\}$. Define $Z'_s$ similarly. Then $Z_s\cong Z'_s$. 
Let
\begin{equation}\label{eq:mr}
M_s = tI-f(H)\pm e_{u_1}e_{v_1}^T\pm \cdots \pm e_{u_{s}}e_{v_{s}}^T, \; s=1,\ldots, r. 
\end{equation} 
Then $e_{u_{r}}^T 
\Bigl(tI-f(A)\pm e_{u_1}e_{v_1}^T \pm\cdots \pm
e_{u_{r-1}}e_{v_{r-1}}^T\Bigl)^{-1}
e_{v_{r}}$ can be written as 
$e_{u_{r}}^T
M_{r-1}^{-1}
e_{v_{r}}$. 
Define $M'_s$ similarly. 
Assume the result holds for $r=k$, that is,   
 $e_{u_{k}}^T
M_{k-1}^{-1}
e_{v_{k}} = e_{u'_{k}}^T
(M')_{k-1}^{-1}
e_{v'_{k}}$. 
Then   
\begin{align*}
&e_{u_{k+1}}^T 
M_k^{-1}
e_{v_{k+1}}&\\
=&e_{u_{k+1}}^T 
\Big(M_{k-1}\pm e_{u_{k}}e_{v_{k}}^T\Big)^{-1}
e_{v_{k+1}}&\\
=&e_{u_{k+1}}^T  
\bigg(M_{k-1}^{-1} \mp \frac{M_{k-1}^{-1}e_{u_{k}}e_{v_{k}}^T 
M_{k-1}^{-1}}
{1+e_{v_{k}}^T  M_{k-1}^{-1}e_{u_{k}}}\bigg)
e_{v_{k+1}}  \quad (\text{by } \eqref{eq:rank1inver})&\\
=&e_{u_{k+1}}^T  M_{k-1}^{-1} e_{v_{k+1}}
 \mp \frac{
 \Big(e_{u_{k+1}} M_{k-1}^{-1} e_{u_{k}}\Big)
\Big(e_{v_{k}}^T  M_{k-1}^{-1}e_{v_{k+1}} \Big) }
{1+e_{v_{k}}^T  M_{k-1}^{-1}e_{u_{k}}}&\\
=&e_{u'_{k+1}}^T  (M'_{k-1})^{-1} e_{v'_{k+1}}
 \mp \frac{
 \Big(e_{u'_{k+1}} (M'_{k-1})^{-1} e_{u'_{k}}\Big)
\Big(e_{v'_{k}}^T  (M'_{k-1})^{-1}e_{v'_{k+1}} \Big) }
{1+e_{v'_{k}}^T  (M'_{k-1})^{-1}e_{u'_{k}}}&\\
=&e_{u'_{k+1}}^T 
(M'_k)^{-1}
e_{v'_{k+1}},&
\end{align*}
where the second last equality follows from the induction hypothesis. Making use of Lemma~\ref{1walkregcomple}, the result about $\bar H$ and $\bar{H'}$ can be proved similarly.\qed

Note that by taking $X'=X$ in the Lemma, we see that for different ways to embed the vertex set of a small (multi)graph $Z$ into a clique of $X$, the value of $e_{u_{r}}^T 
\Big(tI-f(H)\pm e_{u_1}e_{v_1}^T \pm\cdots \pm
e_{u_{r-1}}e_{v_{r-1}}^T\Big)^{-1}
e_{v_{r}}$ is the same. 

Recall that for a subgraph $Y$ of $X$, we use $X\backslash E(Y)$ to denote the graph obtained from $X$ by deleting the edges of $Y$.  
\begin{theorem}\label{thm:cosALQsubg}
Let $X_1$ and $X_2$ be two cospectral 1-walk regular graphs. 
For $i=1,2$,  let $Y_i$ be a graph contained in a clique of $X_i$. Assume $Y_1\cong Y_2$.  
Then $X_1\backslash E(Y_1)$ and $X_2\backslash E(Y_2)$, and their complements, 
are cospectral
with respect to the adjacency, the Laplacian,  the unsigned Laplacian, and the normalized Laplacian matrices. \\
Furthermore, $X_1\backslash E(Y_1)$ and $X_2\backslash E(Y_2)$ are not isomorphic if $X_1$ are $X_2$ are not. 
\end{theorem}

\proof 
For $i=1,2$, let $\hat{Y}_i$ be the graph 
obtained from $Y_i$ by adding $|V(X)|-|V(Y_i)|$ isolated vertices 
so that $V(\hat Y_i)=V(X)$. 
Assume $Y_i$  has $m$ edges with   
 $E(\hat{Y}_1)=\{a_1b_1\ldots, a_mb_m\}$  and 
 $E(\hat{Y}_2)=\{a'_1b'_1,\ldots, a'_mb'_m\}$.  
 Denote the Laplacian matrix of $X_i$ by $L_i$. 
We prove the result for Laplacian matrix. 
That is 
\begin{align*}
&\det\left(tI-L_1+(e_{a_1}-e_{b_1})(e_{a_1}-e_{b_1})^T+\cdots +(e_{a_m}-e_{b_m})(e_{a_m}-e_{b_m})^T\right) \\
= &\det\left(tI-L_2+(e_{a'_1}-e_{b'_1})(e_{a'_1}-e_{b'_1})^T+\cdots +(e_{a'_m}-e_{b'_m})(e_{a'_m}-e_{b'_m})^T\right). 
\end{align*}
We prove this by induction. 
Assume the map taking $a_i$ to $a'_i$ and $b_i$ to $b'_i$ for $i=1,\ldots, m$ is an isomorphism between $\hat Y_1$ and $\hat Y_2$. 
When $m=1$, 
\begin{align*}
&\det\big(tI-L_1+(e_{a_1}-e_{b_1})(e_{a_1}-e_{b_1})^T\big)\\
=&\det(tI-L_1)\det\left(I+(tI-L_1)^{-1}(e_{a_1}-e_{b_1})(e_{a_1}-e_{b_1})^T\right)\\
=&\det(tI-L_1) 
\left(1+(e_{a_1}-e_{b_1})^T(tI-L_1)^{-1}(e_{a_1}-e_{b_1})\right) \quad (\text{ by }\eqref{eq:i+uv^T})\\
=&\det(tI-L_1) 
\Big(1+e_{a_1}^T(tI-L_1)^{-1}e_{a_1} -
e_{a_1}^T(tI-L_1)^{-1}e_{b_1}
-e_{b_1}^T(tI-L_1)^{-1}e_{a_1} \\
& \quad 
+e_{b_1}^T(tI-L_1)^{-1}e_{b_1}\Big)\\
=&\det(tI-L_2) 
\Big(1+e_{a'_1}^T(tI-L_2)^{-1}e_{a'_1} -
e_{a'_1}^T(tI-L_2)^{-1}e_{b'_1}
-e_{b'_1}^T(tI-L_2)^{-1}e_{a'_1} \\
& \quad 
+e_{b'_1}^T(tI-L_2)^{-1}e_{b'_1}\Big)\\ 
=&\det\big(tI-L_2+(e_{a'_1}-e_{b'_1})(e_{a'_1}-e_{b'_1})^T\big),
\end{align*}
where the second last equality follows from the fact $X_1$ and $X_2$ are (Laplacian) cospectral and  Lemma~\ref{lem:valueinclique}.

For $r=1,\ldots, m$, let $Y_{1,r}$ be the graph formed by $a_1b_1, \ldots, a_rb_r$. Define $Y_{2,r}$ similarly. Then $Y_{1,r}\cong Y_{2,r}$.
For $i=1,2$, define $M_{i,s}$ as in \eqref{eq:mr}. 
With $H=L$, $f(x)=x$, and some proper choice of $\pm$ signs and choice of vertices $u_i$ and $v_i$ in \eqref{eq:mr}, we have 
\[
M_{1,4r} = tI-L_1+(e_{a_1}-e_{b_1})(e_{a_1}-e_{b_1})^T+\cdots +(e_{a_r}-e_{b_r})(e_{a_r}-e_{b_r})^T.  
\]

Now the Laplacian characteristic polynomial of the graph obtained from $X_1$ by deleting the edges of $\hat{Y}_1$ satisfy  
\begin{align*}
&\det\left(tI-L_1+(e_{a_1}-e_{b_1})(e_{a_1}-e_{b_1})^T+\cdots +(e_{a_m}-e_{b_m})(e_{a_m}-e_{b_m})^T\right)\\
& =\det\big(M_{1,4(m-1)}+(e_{a_m}-e_{b_m})(e_{a_m}-e_{b_m})^T\big) \\
&=\det(M_{1,4(m-1)})
\bigg(1+ (e_{a_m}-e_{b_m})^T M_{1,4(m-1)}^{-1} (e_{a_m}-e_{b_m})\bigg) \quad (\text{ by }\eqref{eq:i+uv^T})\\
&=\det(M_{1,4(m-1)})
\bigg(1+ e_{a_m}^TM_{1,4(m-1)}^{-1}e_{a_m} - 
e_{a_m}^TM_{1,4(m-1)}^{-1}e_{b_m}
- e_{b_m}^TM_{1,4(m-1)}^{-1}e_{a_m} \\ 
&\quad + e_{b_m}^TM_{1,4(m-1)}^{-1}e_{b_m}
\bigg)\\
&=\det(M_{2,4(m-1)})
\bigg(1+ e_{a'_m}^TM_{2,4(m-1)}^{-1}e_{a'_m} - 
e_{a'_m}^TM_{2,4(m-1)}^{-1}e_{b'_m}
- e_{b'_m}^TM_{2,4(m-1)}^{-1}e_{a'_m} \\ 
&\quad + e_{b'_m}^TM_{2,4(m-1)}^{-1}e_{b'_m}
\bigg)\\
&=\det\left(tI-L_2+(e_{a'_1}-e_{b'_1})(e_{a'_1}-e_{b'_1})^T+\cdots +(e_{a'_m}-e_{b'_m})(e_{a'_m}-e_{b'_m})^T\right), 
\end{align*}
where in the second last equality, $\det(M_{1,4(m-1)})=\det(M_{2,4(m-1)})$ follows from induction hypothesis and the equality of the second factor follows from Lemma~\ref{lem:valueinclique}. 
The result for Laplacian matrix follows. 
 
As above, with proper choices of signs and vertices in \eqref{eq:mr}, 
we can prove the case for adjacency matrix and unsigned Laplacian matrix similarly, 
where deleting an edge corresponds to 
$tI-A+e_ae_b^T+e_be_a^T$ in the adjacency case, 
and corresponds to $tI-S+(e_a+e_b)(e_a+e_b)^T$ in the unsigned Laplacian case.

For the normalized Laplacian case, 
even though the perturbation of the matrix $N=I -\frac1d A$ due to the removal of an edge $ab$ is not simply a linear combination of the rank one matrices $e_ue_v^T$ for $u,v\in \{a,b\}$,  
 the proof argument works similarly. 
Let $X$ be a connected graph with no vertices of degree 1. 
Since the normalized Laplacian $N(X)=D^{-\frac12}L(X)D^{-\frac12}$ and $D^{-1}L(X)$ are similar: 
\[ 
D^{-1}L(X)=D^{-\frac12}N(X)D^{\frac12},  
\]
and 
$D^{-1}L=D^{-1}(D-A)=I-D^{-1}A$, 
we can prove the claim by showing that the characteristic polynomial of 
$D(X\backslash e)^{-1}A(X\backslash e)$ are all the same no matter which edge $e=ab$ we remove from $X$. Note that $D(X)=dI$ here. 
Now 
\begin{align}
&\det \left( tI - D(X\backslash e)^{-1}A(X\backslash e) \right)\nonumber \\
 =& \det \Bigg[ tI -  D(X)^{-1}A(X) + 
e_a \left((\frac1d-\frac{1}{d-1})Ae_a + \frac{1}{d-1} e_b\right)^T \nonumber\\
& \quad +e_b \left((\frac1d-\frac{1}{d-1})Ae_b + \frac{1}{d-1} e_a\right)^T
 \Bigg]  \nonumber \\
&= \det\big[E+e_b \left((\frac1d-\frac{1}{d-1})Ae_b + \frac{1}{d-1} e_a\right)^T\big] \nonumber\\
=& \det E \Big[1+ \big((\frac1d-\frac{1}{d-1})Ae_b+ \frac{1}{d-1} e_a\big)^T E^{-1} e_b\Big], \label{eq:norLper}
\end{align}
where $E=tI -  D(X)^{-1}A(X) + 
e_a \left((\frac1d-\frac{1}{d-1})Ae_a + \frac{1}{d-1} e_b\right)^T$. 
Calculating $\det E$ as above and $E^{-1}$ using Theorem~\ref{thm:shermanmorrison}, we have 
\[
\det E=\det(tI-d^{-1}A)\Big[1+
\left((\frac1d-\frac{1}{d-1})Ae_a + \frac{1}{d-1} e_b\right)^T 
(tI-d^{-1}A)^{-1} e_a\Big]
\]
and 
\[
E^{-1}=(tI-d^{-1}A)^{-1}-
\frac{(tI-d^{-1}A)^{-1}e_a \left((\frac1d-\frac{1}{d-1})Ae_a + \frac{1}{d-1} e_b\right)^T(tI-d^{-1}A)^{-1}}{1+\left((\frac1d-\frac{1}{d-1})Ae_a + \frac{1}{d-1} e_b\right)^T(tI-d^{-1}A)^{-1}e_a}.
\]
Substituting the expressions of $\det E$ and $E^{-1}$ into \eqref{eq:norLper} and simplify, we can express $\det \left( tI - D(X\backslash e)^{-1}A(X\backslash e) \right)$ in terms of $e_u^T (tI-d^{-1}A)^{-1} e_v$ and $e_u^T A(tI-d^{-1}A)^{-1} e_v$ for $u,v\in \{a,b\}$. Now the proof follows similarly as for the Laplacian case and by applying Lemma~\ref{lem:valueinclique} to the functions $f(A)=(tI-d^{-1}A)^{-1}$ and $g(A)=A(tI-d^{-1}A)^{-1}$. 

The non-isomorphism result follows from the fact that an isomorphism between $X_1\backslash E(Y_1)$ and $X_2\backslash E(Y_2)$ is also an isomorphism between $X_1$ and $X_2$. \qed

\begin{remark}
{\normalfont
Let $X_1$ and $X_2$ be two cospectral 1-walk regular graphs on $n$ vertices. For $i=1,2$, let $H_i$ be the adjacency, Laplacian, unsigned Laplacian, or normalized Laplacian matrix of $X_i$. 
Let $f(x)$ be a function defined on the eigenvalues of $H_i$. 
Let $B$ be a block-diagonal matrix of size $n$ such that its unique nonzero block corresponds 
 to a clique of $X_1$. 
By use of Lemma~\ref{lem:valueinclique}, 
the argument in Theorem~\ref{thm:cosALQsubg}
can in fact be used to prove that 
if $P$ is a permutation matrix such that all the non-zero entries of 
$PBP^T$ correspond to a clique of $X_2$, then 
$\det\left(tI-f(H_1)-B\right)=\det\left(tI-f(H_2)-PBP^T\right)$.

 Let $\mathcal{A}=\{A_0=I,A_1,\ldots, A_d\}$ be a symmetric association scheme, with the associated graphs denoted by $G_i$ for $i=1,\ldots, d$. 
We can similarly show that different ways of 
adding an edge of a graph $G_k$ to the graph $G_j$ 
results in two graphs cospectral with respect to the four matrices. 
The result holds in general when adding edges of a small graph contained in a clique of $G_k$ to $G_j$. 
}
\end{remark}

Given two adjacency cospectral graphs $X$ and $X'$, two subsets $U\subseteq V(X)$ and $U'\subseteq V(X')$ are removal-cospectral if  there is a one-to-one mapping $U\rightarrow U'$ such that for any $W\subset U$, the graphs $X\backslash W$ and $X' \backslash W'$ are adjacency cospectral, where $W'$ is the image of $W$ under this mapping \cite{DvF2012}.  
Let $X$ and $X'$ be two cospectral 1-walk regular graphs. 
Let $U$ be a clique of $X$ and $U'$ a clique of $X'$. 
Combining Lemma~\ref{lem:cosp1walkreg}, Theorem~\ref{thm:1walkrechar}, and  
equation~\ref{eq:Jacobithm}, we know that if $U$ and $U'$ are of the same size,  then they are removal-cospectral. Any bijection between $U$ and $U'$ can serve as the mapping required by removal-cospectrality. 
Dalf\'o, van Dam and Fiol showed that certain perturbations (for example, edge deletion) in removal-cospectral subsets of almost distance regular graphs result in adjacency cospectral graphs, and the corresponding subsets in the perturbed graphs are still
removal-cospectral \cite{DvF2012}.  
Hence, restricted to the adjacency case, our above result about different ways of removing edges of a subgraph from a clique of a 1-walk regular graph resulting in adjacency cospectral graphs was an implicit result in \cite{DvF2012}, as one can remove one edge (and its image under the mapping for removal-cospectrality) at one step while keeping the two updated graphs adjacency cospectral and keeping the two subsets removal-cospectral.

\section{Some examples} \label{sec:sgr25,12,5,6}
There are exactly 15 non-isomorphic strongly regular graphs with parameters 
SRG(25,12, 5, 6). 
Their adjacency matrices can be found at Spence's website: 
\url{http://www.maths.gla.ac.uk/~es/srgraphs.php}. 
In Table~\ref{tab:dataek3p3k4}, 
we denote these graphs as $X_i$ for $i=0, 1, \ldots, 14$, 
in accordance with the order of the adjacency matrices given on the website.   
Exactly one of these graphs, $X_{14}$,  is edge-transitive, 
that is, the graphs obtained by deleting an edge from $X_{14}$ 
are all isomorphic. 
(We remark that $X_{14}$ is the Latin square graph 
corresponding to the addition table of $\mathbb{Z}_5$.)

All the 15 graphs have clique number 5: 
with $X_0,\ldots, X_{12}$ containing exactly 3 cliques of size 5, 
and $X_{13}$ and $X_{14}$ containing
exactly 15 cliques of size 5. 
Denote the complete graph on $n$ vertices as $K_n$ and the path graph on $n$ vertices as $P_n$. 
Let $Y$ be $K_2$, $K_3$, $K_4$, $P_3$, or the two graphs on 5 vertices shown in Table~\ref{tab:dataek3p3k4}. 
In the table, 
for each of the 15 strongly regular graphs $X_i$,  
we provide the number of pairwise non-isomorphic graphs we obtain by removing the edges of $Y$ from a clique of $X_i$ in different ways, as the entry in the row indexed by $X_i$ and the column indexed by $Y$. (The resulting subgraphs  are cospectral with respect to $A$, $L$, $S$, and $N$.) 
For example, 
removing an edge from $X_0$ gives a family of 81 pairwise non-isomorphic graphs. 
Similarly, different ways of removing edges of a triangle from $X_0$ gives a family of 132 non-isomorphic cospectral graphs.

Assume $i\neq j\in \{0,1\ldots, 14\}$. 
Let $Y_i$ be a small graph contained in a clique of $X_i$ and $Y_j$ a small graph contained in a clique of $X_j$. 
If $Y_i\cong Y_j$, 
then by Theorem~\ref{thm:cosALQsubg}, 
 $X_i\backslash E(Y_i)$ and $X_j\backslash E(Y_j)$ are cospectral with respect to all the four matrices and are non-isomorphic. 
We provide in the last row of Table~\ref{tab:dataek3p3k4}, the total number of non-isomorphic (cospectral) graphs obtained from these 15 graphs by deleting edges of the small graph $Y$ from a clique of some $X_i$. 
Each of these numbers is obtained by summing up the 15 numbers for the 15 graphs in the column.

\begin{table}[h]
\centering
\begin{tabular}{|c|c|c|c|c|c|c|}
\hline
\diagbox[width=11.2em]{Graphs}{Subgraphs removed \\ from a clique\\ of $X_i$\\ 
}
       &$K_2$             &$K_3$              &$K_4$             &$P_3$
&
 \begin{tikzpicture}[-,>=stealth',shorten >=0.5pt,auto, semithick, scale=0.3]
\tikzstyle{every state}=[draw,circle,thick, fill=black, minimum size=5pt, inner sep=0pt]
    \node[state] (1)  at (0,0){ };  
    \node[state] (2)  at (0,2){ };  
    \node[state] (3)  at (2,1){ }; 
    \node[state] (4)  at (-2,0){ }; 
     \node[state] (5)  at (-2,2){ };

  \draw(2,1)--(0,0);
  \draw(2,1) --(0,2);
  \draw(0,0)--(0,2);
   \draw(0,0)--(-2,0);
   \draw(0,0)--(-2,2);
    \draw(0,2)--(-2,0);
   \draw(0,2)--(-2,2);
   
\end{tikzpicture} 
&
 \begin{tikzpicture}[-,>=stealth',shorten >=0.5pt,auto, semithick, scale=0.4]
\tikzstyle{every state}=[draw,circle,thick, fill=black, minimum size=5pt, inner sep=0pt]
    \node[state] (1)  at (0,0) { }; 
    \node[state] (2)  at (0,2){ };  
    \node[state] (3)  at (-1,1){ }; 
    \node[state] (4)  at (1,1){ }; 
     \node[state] (5)  at (3,1){ };

  \draw(1,1)--(0,2);
  \draw(1,1) --(0,0);
  \draw(1,1)--(-1,1);
   \draw(1,1)--(3,1);
   \draw(0,0)--(0,2);
    \draw(0,0)--(-1,1);
   \draw(0,2)--(-1,1);

\end{tikzpicture} 
\\
\hline
$X_0$	                     &81		  &132		    &48		     &252 
&17  	&9 \\
\hline		
$X_1$	                     &150		  &250		    &89		     &373 
&30  	&15 \\	
\hline
$X_2$	                     &50		  &86		            &32		     &216 
&10  	&5 \\	
\hline
$X_3$	                     &150		  &250		    &89		     &373 
&30  	&15 \\
\hline	
$X_4$	                     &81		  &134		    &50		     &247 
&17  	&9 \\	
\hline
$X_5$	                     &79		  &132		    &46		     &304 
&16  	& 15 \\
\hline	
$X_6$	                     &29		  &50		            &17		     &148 
&6  	&5 \\	
\hline
$X_7$	                     &31		  &50		            &19		     &155 
&7  	&9 \\	
\hline
$X_8$	                     &50		  &86		            &32		     &232 
&10  	&15 \\
\hline	
$X_9$	                     &33		  &52		            &20		     &97 
&8  	&6 \\
\hline	
$X_{10}$	                     &6		  &10		           &6		     &24 
&2  	&3 \\
\hline	
$X_{11}$		                    &79		  &134		   &48		   &298 
&19  	&14 \\
\hline	
$X_{12}$	                     &29		  &49		           &15		     &101 
&6  	&4 \\	
\hline
$X_{13}$	                     &5		  &10		           &5	            &16 
&5  	&4 \\	
\hline
$X_{14}$	                     &1		  &2		           &1		            &5
&1  	&1 \\
\hline
	Total	               & 	854		  & 	1427		         & 	517		            & 	2841
&	184  	&	129 \\
\hline
\end{tabular}
\caption{The number of pairwise non-isomorphic (cospectral) subgraphs of $X_i$ obtained from different ways of deleting edges of a small graph, $K_2$, $K_3$, $K_4$, $P_3$, ect., respectively, from a clique of $X_i$}
\label{tab:dataek3p3k4}
\end{table}

\section{Matrices} 

In the next section, 
we provide a more general way to construct Laplacian 
or unsigned Laplacian cospectral graphs, in the sense that the small subgraphs removed from a clique of a 1-walk regular graph don't have to be isomorphic. 
In this section, 
we develop or review some tools for that. 

\subsection{Similar matrices}

Let $M_1,M_2$ be two cospectral Hermitian matrices. 
We give a characterization of when their rank one updates 
$M_1-v_1v_1^*$ and $M_2-v_2v_2^*$ are also cospectral. 
There is related work in \cite{BNSrank1m,G.Golubmatrix}. However, our result does not follow from it in a straightforward way, hence we present our own proof here. 

\begin{theorem}\label{thm:coscosJ}
Let $M_1$ and $M_2$ be two similar Hermitian matrices of size $n$.
Let $v_1$ and $v_2$ be two vectors in $\mathbb{C}^n$. 
Then the following are equivalent:
\begin{itemize}
\item[(a)]
 $v_1v_1^*-M_1$ and $v_2v_2^*-M_2$ are also similar,
 \item[(b)]
 $
v_1^* (tI-M_1)^{-1}  v_1 = v_2^* (tI-M_2)^{-1} v_2,
 $
\item[(c)]
there exists a unitary matrix $U$ such that 
\[
UM_1=M_2U, \quad \text{ and } \quad Uv_1=v_2.
\]
\end{itemize}
Furthermore, 
if $M_1$, $M_2$, $v_1$ and $v_2$ are all real, 
then the unitary matrix $U$ in Condition $(c)$ can be chosen to be orthogonal. 
\end{theorem}

\proof 
We first prove that Conditions $(a)$ and $(b)$ are equivalent. 
Let $M_1$ and $M_2$ be two similar Hermitian matrices.  
For $i=1,2$, 
\begin{align}
\det\big( tI-(v_iv_i^*-M_i)\big) &= \det(tI+M_i-v_iv_i^*) \nonumber \\
&= \det\Big( (tI+M_i) \big[I-(tI+M_i)^{-1}v_iv_i^*\big] \Big) \nonumber\\
& =  \det(tI+M_i) \big( 1- v_i^* (tI+M_1)^{-1} v_i\big). \label{eq:rank1invdet}
\end{align}
Since the matrices involved are Hermitian, 
we know that  $v_1v_1^*-M_1$ and $v_2v_2^*-M_2$ are similar 
if and only if 
$\det\big( tI-(v_1v_1^*-M_1)\big)=\det\big( tI-(v_2v_2^*-M_2)\big)$.  
Therefore 
by equation \eqref{eq:rank1invdet}, 
for two similar Hermitian matrices $M_1$ and $M_2$, 
their rank one updates $v_1v_1^*-M_1$ and $v_2v_2^*-M_2$ are similar 
 if and only if
$
v_1^* (tI-M_1)^{-1}  v_1 =  v_2^* (tI-M_2)^{-1}  v_2 
 $. 
 
 We now show that $(b)$ implies $(c)$. 
 Let 
 \[
 M_1=\sum_{r=1}^m \theta_r E_r, \quad  M_2=\sum_{r=1}^m \theta_r F_r
 \]
 be the spectral decomposition of $M_1$ and $M_2$, respectively. 
 Then 
 \[
  v_1^*(tI-M_1)^{-1}  v_1 =
  \sum_{r=1}^m  \frac{ v_1^* E_r v_1 }{t-\theta_r}
 \] 
 and Condition $(b)$ holds if and only if 
 \begin{equation}\label{eq:1er1}
  v_1^* E_r v_1=v_2^* F_r v_2, \forall r.
 \end{equation}
Now we construct a unitary matrix $U$ such that 
$UM_1=M_2U$ and $Uv_1=v_2$. 
Any unitary matrix $U'$ that maps 
an orthonormal basis of each eigenspace of $M_1$ 
to an orthonormal basis of the corresponding eigenspace of $M_2$ satisfies  $U'M_1=M_2U'$. 
In choosing a basis for each eigenspace, 
we can start with any unit vector in the eigenspace.  
Choose the first basis vector in the eigenspace associated to $\theta_r$ to be 
$\frac{1}{\sqrt{v_1^* E_r v_1}}E_r v_1$ if $E_r v_1\neq 0$, 
and if $E_r v_1= 0$ we don't put any restrictions on the orthonormal basis of the eigenspace associated to $\theta_r$. 
Choose an orthonormal basis for each eigenspace of $M_2$ in the same way.  
Then the transition matrix $U$ between the two bases  is unitary and satisfies 
$UM_1=M_2U$. Furthermore, and for any $r$ such that $E_r v_1 \neq 0$, we have
\[
U\left(\frac{1}{\sqrt{v_1^* E_r v_1}}E_r v_1\right)
=\frac{1}{\sqrt{v_2^* F_r v_2}}F_r v_2. 
\]
By \eqref{eq:1er1}, 
we conclude that for all $r$, 
$U(E_r v_1)=F_r v_2$. 
Thus 
\[
U v_1=U\sum_r E_r v_1 =\sum_r F_r v_2 =v_2.
\]

Now we prove $(c)$ implies $(a)$. 
Assume $(c)$ holds, 
then 
\[
U(v_1v_1^*-M_1) = v_2v_1^*-UM_1=v_2v_2^*U-M_2U =(v_2v_2^*-M_2)U, 
\]
therefore $v_1v_1^*-M_1$ and $v_2v_2^*-M_2$ are similar.\qed

Let $X_1$ and $X_2$ be two graphs. 
For $i=1,2$, let $M_i$ be the adjacency, Laplacian, or unsigned Laplacian matrix of  $X_i$.  
If $X_1$ and $X_2$, and their complements, are cospectral with respect to $M$, 
we have the following results.   
 Part $(a)$ is due to Johnson and Newman \cite{JohnsonNewmancocos}. 
Denote the all-ones vector by $\mathbf{1}_n$.  
\begin{corollary} \label{cor:coscos1}
\begin{itemize}
\item[(a)] 
If $X$ and $Y$ are adjacency cospectral graphs with cospectral complements, 
then there is an orthogonal matrix $Q$ such that 
\[ 
Q A(X)Q^T = A(Y), \;
Q A(\bar{X})Q^T = A(\bar{Y}) ,\; Q\mathbf{1}=\mathbf{1}. 
\]
\item[(b)] 
If $X$ and $Y$ are unsigned Laplacian cospectral graphs with cospectral complements, 
then there is an orthogonal matrix $Q$ such that 
\[ 
Q S(X)Q^T = S(Y), \;
Q S(\bar{X})Q^T = S(\bar{Y}),\; Q\mathbf{1}=\mathbf{1}. 
\]
\item[(c)] 
If $X$ and $Y$ are Laplacian cospectral graphs, 
then there is an orthogonal matrix $Q$ such that 
\[ 
Q L(X)Q^T = L(Y), \;
Q L(\bar{X})Q^T = L(\bar{Y}),\; Q\mathbf{1}=\mathbf{1}. 
\]
\end{itemize}
\end{corollary}
\proof 
Assume $X$ and $Y$ have $n$ vertices. 
By assumption, 
$A(X)$ and $A(Y)$ are similar, 
and $J-A(X)$ and $J-A(Y)$ are similar. 
Let $M_1=A(X)$, 
$M_2=A(Y)$, 
and $v_1=v_2=\mathbf{1}_n$. 
By the equivalence of condition $(a)$ and $(c)$ in Theorem~\ref{thm:coscosJ}, 
there exists an orthogonal matrix $Q$ such that 
\[
	QA(X)=A(Y)Q,\; Q\mathbf{1}_n=\mathbf{1}_n. 
\]
Therefore 
\[ 
	QA(\bar{X})Q^T=Q(J-I-A(X))Q^T=J-I-A(Y)=A(\bar{Y}). 
\]
The proof for $(b)$ follows similarly with 
$S(\bar{X})=J + (n-2)I-S(X)$. 
The proof for $(c)$ follows from the fact if two graphs are Laplacian cospectral then so are their complements, $L(\bar{X})=nI - J -L(X)$, and from a similar argument as above. 
\qed

\begin{remark} \label{rm:weightggcq}
{\normalfont
Let $\alpha,\beta$ be fixed real numbers. 
Let $X$ be a graph and consider the weighted matrix $A_{\alpha,\beta}(X)=\alpha D(X)+\beta A(X)$. 
We say two graphs are $A_{\alpha,\beta}$-cospectral if their associated 
$A_{\alpha,\beta}$ matrices have the same characteristic polynomials. 
Then a similar argument as in the above corollary shows that:
if $X$ and $Y$ are $A_{\alpha,\beta}$-cospectral 
with cospectral complements, 
then there is an orthogonal matrix $Q$ such that 
\[ 
Q A_{\alpha,\beta}(X)Q^T = A_{\alpha,\beta}(Y), \;
Q A_{\alpha,\beta}(\bar{X})Q^T = A_{\alpha,\beta}(\bar{Y}),\; Q\mathbf{1}=\mathbf{1}. 
\]
}
\end{remark}

\subsection{Gram matrices}

Given a graph $X$, 
if we assign a direction to each edge we obtain an oriented graph $\tilde{X}$. 
For an arc $(a,b)$, 
we call $a$ its tail and $b$ its head. 
The incidence matrix of an oriented graph $\tilde{X}$ is the $(0,\pm1)$-matrix 
with rows indexed by the vertices and columns indexed by the arcs,
such that the $ae$-entry is equal to 1 if vertex $a$ is the head of the arc $e$, 
$-1$ if $a$ is the tail of $e$, 
and 0 otherwise. 
This incidence matrix is called an \textsl{oriented incidence matrix} of $X$. 

Different orientations of $X$ result in different oriented incidence matrices of $X$, 
but for any oriented incidence matrix $B$ of $X$, 
we have $BB^T=L(X)$. 
Furthermore, 
different oriented incidence matrices of the same graph are related by an orthogonal matrix.

The following result is known. For example, it follows from Theorem 10 in Chapter 8 of \cite{HoffmanKunze} by defining $T$ to be the linear transformation mapping $Be_i$ to $Ce_i$ for $i=1,\ldots, n$. For completeness, we present some proofs here. 
\begin{theorem}\label{thm:congruence}
Let $B$ and $C$ be $m\times n$ matrices. 
Then there is a unitary matrix $Q$ such that $QB=C$ 
if and only if $B^*B=C^*C$. \\
If $B$ and $C$ are real, 
then there is an orthogonal matrix $Q$ such that $QB=C$ 
if and only if $B^TB=C^TC$. \qed
\end{theorem}

\proof \textbf{1, SVD}: 
We prove the result for real matrix case. 
Since $B^TB=C^TC$ is a real symmetric matrix, 
it is orthogonally diagonalizable, 
say by $U=[ u_1\;\ldots \; u_n]$ 
to $\Lambda=\diag(\lambda_1,\ldots, \lambda_n)$ 
with $\lambda_1\geq \lambda_2\geq\cdots \geq \lambda_n$. 
That is, 
$U^T B^TB U = \Lambda$. 
Assume $\rk (B)=r$, 
then $\lambda_1\geq\cdots\geq \lambda_r>0$, 
and $\lambda_{r+1}=\cdots=\lambda_n=0$. 

Let 
\begin{equation} \label{eq:lrsvec}
v_i=\frac {1}{\sqrt{\lambda_i}} B u_i, \quad i=1,\ldots, r
\end{equation}
and let 
\[ v_{r+1},\, \ldots, \, v_m
\]
be an orthonormal basis of the null space of $B^T$. 
Then 
\[
V=[v_1\;\ldots \; v_m]
\]
is an orthogonal matrix, and 
\begin{align*}
B U &= [Bu_1\; \cdots \; Bu_r \; Bu_{r+1}\; \cdots\; Bu_m]\\
& =[\sqrt{\lambda_1} v_1\; \cdots \sqrt{\lambda_r} v_r\; 0\; \cdots \; 0] \quad 
\text{ by \eqref{eq:lrsvec}}\\
&= V \Sigma,
\end{align*}
where $\Sigma$ is the $m\times n$ matrix with $\Sigma_{i,i}=\sqrt{\lambda_i}$ 
for $i=1,\ldots, r$ and zero elsewhere. 
That is, 
\begin{equation} \label{eq:Bsvd}
B= V \Sigma U^T
\end{equation}
is a singular value decomposition of $B$. 
Similarly, 
let $w_i =\frac {1}{\sqrt{\lambda_i}} C u_i, \, i=1,\ldots, r$, 
and let $w_{r+1},\ldots, w_m$ be an orthonormal basis of the null space of $C^T$. 
Then for  $W=[w_1,\, \ldots, \, w_m]$, 
\begin{equation}\label{eq:Csvd}
C=W\Sigma U^T
\end{equation}
 is a singular value decomposition of $C$. 
From \eqref{eq:Bsvd} and \eqref{eq:Csvd} we have 
\[ 
C=W\Sigma U^T =W(V^TBU)U^T=(WV^T)B. 
\]
Therefore $Q=WV^T$ is an orthogonal matrix that satisfies $QB=C$.\qed

\proof \textbf{2, reflection induction}: We prove the result for real matrix case. 
Let the columns of $B$ and $C$ be respectively $\seq b1n$ and $\seq c1n$. Assume that
$\rk(B)=r$ and that $\seq b1r$ is a basis for the column space of $B$. Then $\seq c1r$
is a basis for the column space of $C$. 

Since $b_1$ and $c_1$ have the same length, the matrix $Q_1$ representing reflection in
the hyperplane $(b_1-c_1)^\perp$ is an orthogonal matrix swapping $b_1$ and $c_1$ and
\[
	(Q_1B)^TQ_1B = B^TB = C^TC; 
\]
we obtain a matrix $Q_1B$ that share the same first column as $C$
and is equivalent to $C$ (since $(Q_1B)^TQ_1B  = C^TC$). We denote $Q_1B$ as $B$. 

Now assume inductively that $b_i=c_i$ for $i=1,\ldots,k$, with $1\le k\le r$.
If $y$ and $z$ are two vectors such that $\ip{y}y=\ip{z}z$ and
\[
	\ip{c_i}y = \ip{c_i}z,\quad (i=1,\ldots,k)
\]
then $y-z$ is orthogonal to $\seq c1k$ and the reflection in $(y-z)^\perp$ fixes
$\seq c1k$ and swaps $y$ and $z$. If $k<r$, 
take $y=Be_{r+1}$ and $z=Ce_{r+1}$,
and the above implies that 
there is an orthogonal
matrix $Q_{k+1}$ such that the first $k+1$ columns of $Q_{k+1}B$ and $C$ are equal.

To complete the proof, we observe that if the first $r$ columns of $B$ is a basis of $\col(B)$
and are equal to the first $r$ columns of $C$, 
then $B^TB=C^TC$ implies $B=C$. The theorem follows.\qed

\section{More (unsigned) Laplacian cospectral graphs}

In Section~\ref{sec:samesubgraph}, 
we showed that for cospectral 1-walk regular graphs $X$ and $X'$, the graph obtained by removing the edges of a small graph $Y$ from a clique of $X$ and the graph obtained by removing the edges of $Y'\cong Y$ from a clique of  $X'$ are  cospectral (with cospectral complements) with respect to $A,L,S$ and $N$.    
In fact, 
there is more we can say about the Laplacian  and unsigned Laplacian case. 
As shown below, the isomorphism condition on the removed subgraphs $Y$ and $Y'$ can be weakened to certain cospectrality.  
Recall that for a subgraph $Y$ of $X$, we use $X\backslash E(Y)$ to denote the graph obtained from $X$ upon removing the edges of $Y$. 

\begin{theorem}\label{thm:cosLorQcossubg}
Let $X_1$ and $X_2$ be two cospectral 1-walk regular graphs. 
For $i=1,2$, let $Y_i$ be a graph contained in a clique of $X_i$. 
\begin{itemize}
\item[(a)]
If $Y_1$ and $Y_2$ are Laplacian cospectral, 
then so are $X_1\backslash E(Y_1)$ and $X_2\backslash E(Y_2)$. 
\item[(b)]
If $Y_1$ and $Y_2$, and their complements, are unsigned Laplacian cospectral, 
then so are $X_1\backslash E(Y_1)$ and $X_2\backslash E(Y_2)$ (their complements are cospectral as well).  
\end{itemize}
\end{theorem}

\proof 
(a) Assume $|V(X_1)|=|V(X_2)|=n$ and $|V(Y_1)|=|V(Y_2)|=m$. 
Reorder the vertices of $X_1$ and $X_2$ if necessary, 
we can assume without loss of generality that for $i=1,2$, the clique of $X_i$ that contains $Y_i$ consists of vertices $1,\ldots,m$.

For $i=1, 2$, 
let $B_i$ be a signed incidence matrix of $Y_i$, 
that is, $B_iB_i^T=L(Y_i)$. 
Assume $Y_1$ and $Y_2$ are Laplacian cospectral. 
By Corollary~\ref{cor:coscos1} $(c)$, 
there exist an orthogonal matrix $Q$ such that 
$Q^T L(Y_1)Q=L(Y_2)$ and $Q\mathbf{1}=\mathbf{1}$. 
Hence $Q^TB_1B_1^TQ=B_2B_2^T$. 
By Theorem~\ref{thm:congruence}, 
there exists an orthogonal matrix $Q_0$
such that 
\begin{equation}\label{eq:QBQC}
	B_2=Q^TB_1Q_0.
\end{equation}

For $i=1,2$, let $\hat{Y}_i$ be the graph obtained from $Y_i$ by adding $n-m$ isolated vertices so that $\hat{Y_i}$ has the same vertex set as $X$, and let 
\[
\hat{B}_i=\begin{bmatrix} B_i\\ 0\end{bmatrix}
\]
be the matrix obtained from $B_i$ by adding $n-m$ rows of zero. 
That is, $\hat B_i$ is a signed incidence matrix of $\hat Y_i$. 
Then $\hat{B}_i\hat{B}_i^T =L(\hat{Y}_i)$. 
Our goal is to prove that 
the two graphs obtained by removing the edges of $\hat{Y}_1$ from $X_1$ and by removing the edges of $\hat{Y}_2$ from $X_2$ are Laplacian cospectral. 
That is, 
\[
\det \big(tI-L(X_1)+L(\hat{Y_1})\big) = 
\det \big(tI-L(X_2)+L(\hat{Y_2})\big). 
\]
Since 
\begin{align*}
\det \big(tI-L(X_1)+L(\hat{Y_1})\big)  
&= \det \Big( \big(tI-L(X_1)\big) \big(I+\left(tI-L(X_1)\right)^{-1} \hat{B}_1\hat{B}_1^T\big)\Big)\\
&= \phi_L(X_1,t) \det \big(I +\hat{B}_1^T (tI-L(X_1)\big)^{-1} \hat{B}_1\big) \;\; (\text{ by Lemma }\ref{lem:i-cddcddet}) 
\end{align*}
and $\phi_L(X_1,t)=\phi_L(X_2,t)$, it is equivalent to prove 
\[
\det\big(I+\hat{B}_1^T (tI-L(X_1)\big)^{-1} \hat{B}_1\big)
=\det\big(I+\hat{B}_2^T (tI-L(X_2)\big)^{-1} \hat{B}_2\big).
\]
Since $X_1$ and $X_2$ are cospectral 1-walk regular graphs and the first $m$ vertices form a clique in both graphs, 
by Lemma~\ref{lem:1walkinver}, 
there exist scalars $\alpha, \beta$ such that for $i=1,2$, 
$(tI-L(X_i))^{-1}$ is of the form 
\[ (tI-L(X_i))^{-1}=\begin{bmatrix}
\alpha I - \beta(J-I) & M_{i,1}\\
M_{i,2} & M_{i,3} \end{bmatrix}, 
\]
for some matrices $M_{i,1},M_{i,2}$ and $M_{i,3}$, 
whose value does not matter here. 
Now 
\begin{align*}
\det\big(I+\hat{B}_2^T (tI-L(X_2)\big)^{-1} \hat{B}_2\big)  
&= \det \Big( I+
\begin{bmatrix} B_2^T & 0\end{bmatrix} 
\begin{bmatrix}
\alpha I - \beta(J-I) & M_{2,1}\\
M_{2,2} & M_{2,3} \end{bmatrix}
\begin{bmatrix} B_2 \\ 0\end{bmatrix}
\Big)\\
&= \det \Big( I+
B_2^T \big(\alpha I - \beta(J-I)\big)
B_2
\Big)\\
&= \det \Big( I+
Q_0^TB_1^TQ \big(\alpha I - \beta(J-I)\big)
Q^TB_1Q_0
\Big) \quad \text{ by \eqref{eq:QBQC}}\\
&= \det \Big(I+ 
B_1^TQ \big(\alpha I - \beta(J-I)\big)
Q^TB_1
\Big)\\
&= \det \Big( I+
B_1^T \big(\alpha I - \beta(J-I)\big)
B_1
\Big) \quad (\text{ since } Q\mathbf{1}=\mathbf{1})\\
&=\det\big(I+\hat{B}_1^T (tI-L(X_1)\big)^{-1} \hat{B}_1\big).  
\end{align*}

(b) By use of the vertex-edge incidence matrix instead of an oriented incidence matrix of $X$, 
and Corollary~\ref{cor:coscos1} (b), 
the result for unsigned Laplacian case follows similarly.\qed

As in Remark~\ref{rm:weightggcq}, 
we mention a more general result on the weighted matrix 
$A_{\alpha,\beta}(X)=\alpha D(X)+\beta A(X)$ with $\alpha\geq |\beta|$.  
In this case, $A_{\alpha,\beta}(X)$ is positive-semidefinite, 
and hence there exists a matrix $B$ such that $BB^T=A_{\alpha,\beta}$. 
With a similar proof as in Theorem~\ref{thm:cosLorQcossubg} we have the following result. 
Recall two graphs $X$ and $Y$ are $A_{\alpha,\beta}$-cospectral if 
$\det\left(tI-A_{\alpha,\beta}(X)\right)=\det\left(tI-A_{\alpha,\beta}(Y)\right)$.

\begin{lemma}
Let $\alpha,\beta$ be real numbers with $\alpha\geq |\beta|$. 
Let $X_1$ and $X_2$ be two cospectral 1-walk regular graphs. 
For $i=1,2$ let $Y_i$ be a graph contained in a clique of $X_i$.   
If $Y_1$ and $Y_2$  are $A_{\alpha,\beta}$-cospectral graphs with cospectral complements,  then so are $X_1\backslash E(Y_1)$ and $X_2\backslash E(Y_2)$. 
 \qed
\end{lemma}

\begin{example}
\begin{figure}[h]
     \centering
     \begin{subfigure}[b]{0.44\textwidth}
         \centering
         \includegraphics[scale=0.6]{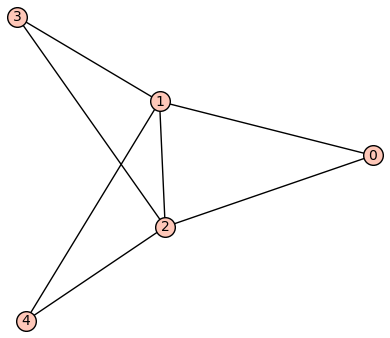}
         \caption{$Y_1$}
         \label{fig:8-1}
     \end{subfigure}
     \hfill
     \begin{subfigure}[b]{0.44\textwidth}
         \centering
         \includegraphics[scale=0.66]{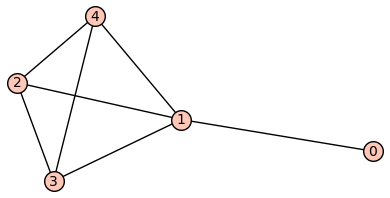}
         \caption{$Y_2$}
         \label{fig:8-2}
     \end{subfigure}
     \caption{A pair of unsigned Laplacian cospectral graphs with cospectral complements}
     \label{fig:unlapcocos}
 \end{figure}  
Let $Y_1$ and $Y_2$ be the two graphs as shown in Figure~\ref{fig:unlapcocos}     
(also in Table~\ref{tab:dataek3p3k4}).
They  are 
unsigned Laplacian cospectral with cospectral complements. 
By Theorem~\ref{thm:cosLorQcossubg}, 
removing edges of $Y_1$ or edges of $Y_2$, respectively,
 from a clique of a 1-walk regular graph gives unsigned Laplacian cospectral graphs, 
which are not isomorphic (they have different degree sequences). 
Therefore, 
for each of the strongly regular graphs $X_i$ in  Table~\ref{tab:dataek3p3k4}, 
we can take the union of the two families of non-isomorphic graphs obtained upon  deleting edges of $Y_1$ or $Y_2$, respectively, from a clique of $X_i$, 
and get a bigger family of pairwise  non-isomorphic but unsigned Laplacian cospectral graphs. 
 We can further take the union of these bigger families of graphs among the 15 strongly regular graphs and obtain an even bigger family of pairwise non-isomorphic but unsigned Laplacian cospectral graphs.
The sizes of the families of non-isomorphic graphs resulting from deleting edges of $Y_1$ or $Y_2$ from a clique of $X$ correspond to the last two columns of Table~\ref{tab:dataek3p3k4}. 
For example, 
deleting edges of $Y_1$ in a clique of $X_0$ gives a family of 17 non-isomorphic but ($A,L,S,N$) cospectral graphs, 
deleting edges of $Y_2$ in a clique of $X_0$ gives such a family of size 9. 
They together give
a family of  $17+9=26$ nonisomorphic unsigned Laplacian cospectral graphs. 
Taking union of these families for the 15 graphs, we have a family of 313 nonisomorphic unsigned Laplacian cospectral graphs.

\end{example}

\begin{example}
Note that for unsigned Laplacian case, 
the condition that the two small graphs are unsigned Laplacian cospectral with cospectral complements is important. 
For example, 
the two graphs $Y_3$ and $Y_4$ in Figure~\ref{fig:unlapcos} are unsigned Laplacian cospectral, 
but don't have cospectral complements. 
Removing their edges inside a clique of   SRG(36,14,7,4) does not always result in unsigned Laplacian cospectral graphs. 
\begin{figure}[h]
     \centering
     \begin{subfigure}[b]{0.44\textwidth}
         \centering
         \includegraphics[scale=0.6]{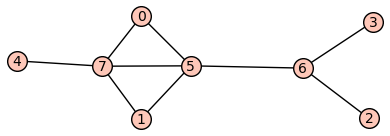}
         \caption{$Y_3$}
         \label{fig:8-1}
     \end{subfigure}
     \hfill
     \begin{subfigure}[b]{0.44\textwidth}
         \centering
         \includegraphics[scale=0.6]{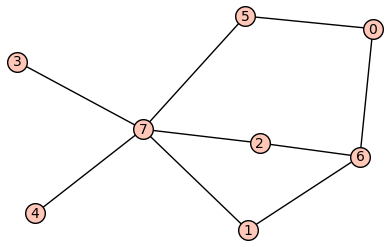}
         \caption{$Y_4$}
         \label{fig:8-2}
     \end{subfigure}
     \caption{A pair of unsigned Laplacian cospectral graphs with non-cospectral complements}
     \label{fig:unlapcos}
 \end{figure}      

\end{example}

\section*{Acknowledgements}
The authors thank Edwin van Dam  for pointing us to some references. 
They also thank the anonymous referees for the very helpful comments and suggestions in improving the exposition of the paper. \\ 
Chris Godsil is supported by Natural Sciences and Engineering Research Council of Canada, Grant No. RGPIN-9439. 
Wangting Sun was a visiting Ph. D. student at the Department of Combinatorics \& Optimization at University of Waterloo from March 2022 to March 2023; supported by China Scholarship Council (No. 202106770031).

\end{document}